\documentclass{amsart}

\newtheorem{theorem}{Theorem}
\newcommand{\bt}{\begin{theorem}}
\newcommand{\et}{\end{theorem}}
\newtheorem{lemma}{Lemma}
\newcommand{\bl}{\begin{lemma}}
\newcommand{\el}{\end{lemma}}
\newtheorem{corollary}{Corollary}
\newcommand{\bc}{\begin{corollary}}
\newcommand{\ec}{\end{corollary}}
\newtheorem{problem}{Problem}
\newcommand{\bprob}{\begin{problem}}
\newcommand{\eprob}{\end{problem}}

\newtheorem{example}{Example}
\newcommand{\bex}{\begin{example}}
\newcommand{\eex}{\end{example}}

\newcommand{\beq}{\begin{equation}}
\newcommand{\eeq}{\end{equation}}
\newcommand{\benum}{\begin{enumerate}}
\newcommand{\eenum}{\end{enumerate}}
\newcommand{\N}{\ensuremath{ \mathbf N }}

\newcommand{\C}{\ensuremath{\mathbf C}}

\newcommand{\mcf}{\ensuremath{ \mathcal F}}
\newcommand{\mcg}{\ensuremath{ \mathcal G}}
\newcommand{\mch}{\ensuremath{ \mathcal H}}
\newcommand{\mci}{\ensuremath{ \mathcal I}}
\newcommand{\mcj}{\ensuremath{ \mathcal J}}

\newcommand{\mbt}{\ensuremath{ \mathbf t}}

\newcommand{\mbw}{\ensuremath{ \mathbf w}}

\newcommand{\mbz}{\ensuremath{ \mathbf z}}

\newcommand{\bsmallmat}{\left(\begin{smallmatrix}}
\newcommand{\esmallmat}{\end{smallmatrix}\right)}

\newcommand{\bmat}{\left(\begin{matrix}}
\newcommand{\emat}{\end{matrix}\right)}

\DeclareMathOperator{\qqand}{\qquad\text{and}\qquad}

\usepackage{fouridx}
\newcommand{\std}[1]{\fourIdx{\circ}{}{}{}{#1}}
\newcommand{\nstd}[1]{\fourIdx{*}{}{}{}{#1}}

\usepackage{amssymb,latexsym}

\title[Continuity of affine varieties]{A nonstandard proof of continuity of affine varieties}
\author{Melvyn B.  Nathanson}
\address{Department of Mathematics\\Lehman College (CUNY)\\Bronx, NY 10468} 
\email{melvyn.nathanson@lehman.cuny.edu}



\date{\today}

\begin{document}

\begin{abstract}
Extending the classical result 
that the roots of a polynomial with coefficients in $\C$ are continuous functions 
of the coefficients of the polynomial, nonstandard analysis is used to prove that  
if $\mcf = \{f_{\lambda}  :\lambda \in \Lambda\}$ 
is a set of polynomials in $\C[t_1,\ldots, t_n]$ and if 
 $\fourIdx{*}{}{}{}{\mcg} = \{g_{\lambda}  :\lambda \in \Lambda\}$ 
 is a set of polynomials in $\fourIdx{*}{}{}{0}{\C} [t_1,\ldots, t_n]$ 
 such that $g_{\lambda}$ is an infinitesimal deformation of $f_{\lambda}$ 
 for all $\lambda \in \Lambda$,  
 then the nonstandard affine variety $\fourIdx{*}{}{}{0}{V}(\mcg)$ 
 is an infinitesimal deformation of the affine variety $V(\mcf)$. 
\end{abstract}

\maketitle

\section{Continuity of the roots of polynomials}

Let  $(a_i)_{i\in I}$ and $(b_i)_{i\in I}$ be sequences of elements 
in the field \C\ of complex numbers and let $\delta > 0$.  
The sequence $(b_i)_{i\in I}$  is a 
\emph{$\delta$-deformation}\index{deformation!$\delta$-deformation} 
of the  sequence $(a_i)_{i\in I}$ if 
\[
|b_i - a_i| < \delta
\]
for all $i \in I$.  

Let $f(t) = \sum_{i=0}^n a_i t^i$ and $g(t) = \sum_{i=0}^n b_i t^i$ be 
polynomials in $\C[t]$ of degree $n$. 
 Let $\delta > 0$.  The polynomial  $g(t)$ is a  
\emph{$\delta$-deformation}\index{deformation} of the polynomial $f(t)$ if 
the sequence $(b_0, b_1,\ldots, b_n)$ of coefficients of $g(t)$ is a 
$\delta$-deformation of the  sequence $(a_0, a_1,\ldots, a_n)$ of coefficients of  $f(t)$.

Every polynomial of degree $n$ in $\C[t]$ 
has $n$ roots, counting multiplicity.  Let 
$(\zeta_1,\ldots, \zeta_{n})$ be the sequence of roots of $f(t)$ 
and let $(\omega_1,\ldots, \omega_n)$ be the sequence of roots of $g(t)$, 
where a root of multiplicity $m$ appears $m$ times in the sequence.  
The polynomial $g(t)$ is \emph{$\varepsilon$-aligned} with the polynomial $f(t)$ 
if there is an ordering  $(\zeta_1,\ldots, \zeta_{n})$ 
of the sequence  of the $n$ roots of $f(t)$ 
and an ordering $(\omega_1,\ldots, \omega_n)$ of the sequence  of the $n$ roots of $g(t)$ 
such that the sequence $(\omega_1,\ldots, \omega_n)$ is an $\varepsilon$-deformation 
of the sequence $(\zeta_1,\ldots, \zeta_{n})$.

A classical theorem on the continuity of roots of a polynomial states that, 
for every $\varepsilon > 0$, 
there exists $\delta > 0$ such that if the polynomial $g(t)$ 
is a $\delta$-deformation of the polynomial $f(t)$, then 
$g(t)$ is \emph{$\varepsilon$-aligned} with $f(t)$.  
Most proofs  of this result require sophisticated mathematical 
tools (cf.~\cite{cuck-corb89}--\cite{hiro20},\cite{mard66},\cite{ostr73},\cite{rahm-schm02}\cite{whit72}).  
Recently, Nathanson and Ross~\cite{nath-ross22} gave a proof 
that used only high school algebra and the definition of continuity.  
Moreover, their proof is valid for polynomials with coefficients 
in every algebraically closed field with an absolute value.  

The theorem on the continuity of the roots of a polynomial 
has an equivalent formulation in nonstandard analysis.   
Let $\fourIdx{*}{}{}{}{\C}$   be the field of nonstandard complex numbers.  We have 
\[
\nstd{\C} =  \fourIdx{*}{}{}{0}{\C} \cup \fourIdx{*}{}{}{\infty}{\C} 
\]
where $\fourIdx{*}{}{}{0}{\C}$   is the ring of finite nonstandard complex numbers 
and $\fourIdx{*}{}{}{\infty}{\C}$ is the set of infinite nonstandard complex numbers. 
Let $\std{w}$ denote the standard part of the nonstandard complex number 
$w \in  \nstd{\C_0}$.  
The function 
$\fourIdx{\circ}{}{}{}{}: \fourIdx{*}{}{}{0}{\C} \rightarrow \C$ is a ring homomorphism 
such that $\fourIdx{\circ}{}{}{}{w} = 0$ if and only if $w$ is infinitesimal.   
We write $w_1 \approx w_2$ if $ \fourIdx{\circ}{}{}{1}{w} =  \fourIdx{\circ}{}{}{2}{w}$ 
or, equivalently, if $ \fourIdx{\circ}{}{}{}{(w_1 - w_2)} = 0$ 
or $w_1 = w_2 + \varepsilon$ for some 
infinitesimal $\varepsilon$. 
We have $ \fourIdx{\circ}{}{}{}{w} = w$ if and only if $w \in \C$. 

Let  $(b_i)_{i\in I}$ and $(c_i)_{i\in I}$ be sequences of elements 
in the field 
 $\fourIdx{*}{}{}{}{\C}$ of nonstandard complex numbers. 
The sequence $(b_i)_{i\in I}$  is an 
\emph{infinitesimal deformation}\index{deformation!infinitesimal-deformation} 
of the  sequence $(c_i)_{i\in I}$  if  for all $i \in I$  there is an infinitesimal $\varepsilon_i$ 
such that 
\[
b_i = c_i + \varepsilon_i. 
\]

The polynomial $g(t) = \sum_{i=0}^n b_i t^i \in \fourIdx{*}{}{}{0}{\C}[t]$ is an 
\emph{infinitesimal deformation} of the polynomial $f(t) = \sum_{i=0}^n a_i t^i \in \C[t]$ 
if the sequence $(b_0, b_1,\ldots, b_n)$ of coefficients of $g(t)$ is an 
infinitesimal deformation of the  sequence $(a_0, a_1,\ldots, a_n)$ of coefficients of  $f(t)$ 
or, equivalently, if $\fourIdx{\circ}{}{}{}{b_i}  = a_i$ 
 for all $i \in \{0,1,\ldots, n\}$. 

The polynomial $g(t)$ is \emph{infinitesimally aligned} with the polynomial $f(t)$ 
if there is an ordering  $(\zeta_1,\ldots, \zeta_{n})$ 
of the sequence  of the $n$ roots of $f(z)$ 
and an ordering $(\omega_1,\ldots, \omega_n)$ of the sequence  of the $n$ roots of $g(t)$ 
such that the sequence $(\omega_1,\ldots, \omega_n)$ is an infinitesimal deformation 
of the sequence $(\zeta_1,\ldots, \zeta_{n})$.  

Ross~\cite{ross22} proved that if the polynomial $g(z)  \in \fourIdx{*}{}{}{0}{\C}[t]$ 
is an infinitesimal deformation 
of the polynomial $f(z) \in \C[t]$, 
then the roots of $g(z)$ are infinitesimally aligned with the roots of $f(z)$.  
This is equivalent to the standard theorem on the continuity of the roots of a polynomial.  

The goal  of this paper is to use nonstandard analysis to extend the theorem 
on the continuity of the roots of a polynomial 
in one variable to the continuity of the roots of a polynomial in several variables 
or of a set of polynomials in several variables, that is, to hypersurfaces and affine algebraic varieties.

\emph{Notation:} The left superscript * is used to distinguish sets in the nonstandard world 
from sets in the standard world.

\section{Hypersurfaces} 
We begin with standard and nonstandard hypersurfaces.  

Let $t_1,\ldots, t_n$ be variables and let $\mbt = (t_1,\ldots, t_n)$.  
Let $\N_0^n$ be the set of $n$-tuples of nonnegative integers.  
For $I = (i_1,\ldots, i_n) \in \N_0^n$, 
we define the monomial 
\[
\mbt^I = t_1^{i_1} \cdots t_n^{i_n}.
\]
Let \mci\ be a finite subset of $\N_0^n$. 
The \emph{hypersurface} defined by the polynomial   
\[
f(\mbt)  = \sum_{I \in \mci} a_I\mbt^I  \in \C[\mbt] 
\]
is 
\[
V(f) = \{\mbz \in \C^n:f(\mbz) = 0\}.
\]
The \emph{nonstandard hypersurface}\index{nonstandard hypersurface}\index{hypersurface!nonstandard} defined by the polynomial 
\[
g(\mbt)  = \sum_{I \in \mci} b_I \mbt^I \in  \fourIdx{*}{}{}{}{\C}[\mbt]  
\] 
is 
\begin{align*}
\fourIdx{*}{}{}{}{V}(g) 
&  = \{\mbw \in \fourIdx{*}{}{n}{}{\C}:  g   (\mbw) \approx 0 \}.
\end{align*}
The \emph{nonstandard finite hypersurface} defined by the polynomial $g(\mbt)$ 
is 
\begin{align*}
\fourIdx{*}{}{}{0}{V}(g) 
&  = \{\mbw \in \fourIdx{*}{}{n}{0}{\C}:  g   (\mbw) \approx 0 \}  \\
& = \{\mbw \in \fourIdx{*}{}{n}{0}{\C}:  \std{g}(\mbw) = 0 \}
\end{align*}
and the   \emph{nonstandard infinite hypersurface} defined by the polynomial $g(\mbt)$ 
is 
\begin{align*}
\fourIdx{*}{}{}{\infty}{V}(g) 
&  = \{\mbw \in  \fourIdx{*}{}{n}{\infty}{\C}  :  g   (\mbw) \approx 0 \}.
\end{align*}
The hypersurfaces $\fourIdx{*}{}{}{0}{V}(g)$ and $\fourIdx{*}{}{}{\infty}{V}(g)$ are disjoint, 
but $\fourIdx{*}{}{}{0}{V}(g) \cup \fourIdx{*}{}{}{\infty}{V}(g)$ may be a proper subset of 
$\fourIdx{*}{}{}{}{V}(g)$.   For example, if $g(t_1,t_2) = t_1t_2-1$, then, for every nonzero 
infinitesimal $\varepsilon$, we have  $(\varepsilon,1/\varepsilon) \in \fourIdx{*}{}{}{}{V}(g)$ 
but $(\varepsilon,1/\varepsilon) \notin \fourIdx{*}{}{}{0}{V}(g) \cup \fourIdx{*}{}{}{\infty}{V}(g)$.

Consider polynomials 
\[
g = g(\mbt) = \sum_{I \in \mci} b_I \mbt^I \in  \fourIdx{*}{}{}{}{\C}[\mbt] 
\]
and 
\[
h = h(\mbt) = \sum_{I \in \mci} c_I \mbt^I \in  \fourIdx{*}{}{}{}{\C}[\mbt]. 
\]
The polynomial $g(\mbt)$ is an infinitesimal deformation 
of the polynomial $h(\mbt)$, 
denoted $g \approx h$,  
if the sequence $(b_I)_{I \in \mci}$ is an infinitesimal deformation of the sequence 
$(c_I)_{I\in \mci}$, that is, if there is a sequence of infinitesimals $(\varepsilon_I)_{I \in \mci}$ 
such that 
\[
b_I = c_I + \varepsilon_I
\]
for all $I \in \mci$. 
Let $g_1,g_2,h_1,h_2 \in  \fourIdx{*}{}{}{0}{\C}[\mbt]$. 
If $g_1 \approx h_1$ and $g_2 \approx h_2$, then 
\[
g_1+g_2 \approx h_1 + h_2 \qqand g_1g_2 \approx h_1 h_2. 
\]

For polynomials with finite coefficients $g(\mbt)  \in  \fourIdx{*}{}{}{0}{\C}[\mbt]$ 
and $h(\mbt)  \in  \fourIdx{*}{}{}{0}{\C}[\mbt]$, 
we define 
\[
 \fourIdx{\circ}{}{}{}{g}(\mbt) =  \sum_{I \in \mci} \std{b_I} \mbt^I \in \C[\mbt] 
 \qqand 
  \fourIdx{\circ}{}{}{}{h}(\mbt) =  \sum_{I \in \mci} \std{c_I} \mbt^I \in \C[\mbt].
\]
It follows that $g \approx h$ if and only if $ \fourIdx{\circ}{}{}{}{g} = \fourIdx{\circ}{}{}{}{h} $.
For the polynomial 
\[
f(\mbt) = \sum_{I \in \mci} a_I\mbt^I  \in \C[\mbt] 
\]
we have $g \approx f$ if and only if $ \fourIdx{\circ}{}{}{}{g}(\mbt) = f(\mbt)$ if and only if 
\[
\fourIdx{\circ}{}{}{i}{b} = a_i
\]
for all $i \in \mci$.

The identities 
\[
 \fourIdx{\circ}{}{}{}{(g +h)}(\mbt) = \fourIdx{\circ}{}{}{}{g}(\mbt)  + \fourIdx{\circ}{}{}{}{h}(\mbt)  
\]
and 
\[
 \fourIdx{\circ}{}{}{}{(g h)}(\mbt) = \fourIdx{\circ}{}{}{}{g}(\mbt)   \fourIdx{\circ}{}{}{}{h}  (\mbt)
\]
imply:  

\bl 
\[
^{\circ}:   \fourIdx{*}{}{}{0}{\C}[\mbt] \rightarrow \C[\mbt] 
\]
is a ring homomorphism. 
\el

For $\mbw = (w_1,\ldots, w_n) \in \fourIdx{*}{}{n}{0}{\C}$, we define 
\[
\std{\mbw} = \fourIdx{\circ}{}{}{}{\left(w_1,   \ldots, w_n \right)}    
=  \left( \fourIdx{\circ}{}{}{}{w_1} ,  \ldots, \fourIdx{\circ}{}{}{}{w_n}  \right) .  
\]
If $I = (i_1,\ldots, i_n) \in \N_0^n$, then   
\[
\mbw^I = w_1^{i_1} \cdots w_n^{i_n} 
\]
and 
\[
\std{\left( \mbw^I \right)} = \fourIdx{\circ}{}{}{}{\left(w_1^{i_1} \cdots w_n^{i_n}\right)}    
=  ( \fourIdx{\circ}{}{}{}{w_1} )^{i_1} \cdots ( \fourIdx{\circ}{}{}{}{w_n} )^{i_n} 
= \left(  \std{\mbw}\right)^I. 
\]
Thus, the notation $\std{\mbw^I} = \std{\left( \mbw^I \right)} = \left(  \std{\mbw}\right)^I$ 
is unambiguous. 

For $\mbw = (w_1,\ldots, w_n) \in \fourIdx{*}{}{n}{0}{\C}$ 
and $\mbz = (z_1,\ldots, z_n) \in \C^n$, we have $\mbw \approx \mbz$ if 
$\std{w_i} = z_i$ for all $i \in \{1,\ldots, n\}$. 
Let  $W$ be a subset of $\fourIdx{*}{}{n}{0}{\C}$ and 
let $Z$ be a subset of $\C^n$.  
The set $W$ is an infinitesimal deformation of the set $Z$ 
if 
\[
W = \{ \mbw\in \fourIdx{*}{}{n}{0}{\C} : \mbw \approx \mbz \text{ for some } \mbz \in Z \} 
= \{ \mbw\in \fourIdx{*}{}{n}{0}{\C} : \std{\mbw} \in Z \}. 
\]

\bt                    \label{NonstandardAffine:theorem:hypersurface} 
Let $f(\mbt) \in \C[\mbt]$ be a nonzero polynomial.  
If the polynomial $g(\mbt) \in \fourIdx{*}{}{}{0}{\C}[\mbt]$ is an infinitesimal deformation of $f(\mbt)$, then 
the nonstandard finite hypersurface $\fourIdx{*}{}{}{0}{V}(g) $ is an infinitesimal 
deformation of the  hypersurface $V(f)$.  Equivalently, 
\[
\fourIdx{*}{}{}{0}{V}(g)  = \left\{ \mbw \in \fourIdx{*}{}{n}{0}{\C} : \std{\mbw} \in V(f)\right\}. 
\] 
\et

\begin{proof} 
Let 
\[
f(\mbt) = \sum_{i \in \mci} a_I \mbt^I \in \C[\mbt] 
\]
be a nonzero polynomial and let 
\[
g(\mbt) = \sum_{i \in \mci} b_I \mbt^I \in \fourIdx{*}{}{}{0}{\C}[\mbt]
\]
be an infinitesimal deformation  of $f$.

Let $\mbz  \in V(f)$.  
If $\mbw \in  \fourIdx{*}{}{n}{0}{\C} $ and $\std{\mbw} = \mbz$, 
then 
\begin{align*}
\std{g}(\mbw) 
& = \std{\left( \sum_{I \in \mci} b_I \mbw^I \right)} 
=  \sum_{I \in \mci}      \fourIdx{\circ}{}{}{I}{b}    \fourIdx{\circ}{}{I}{}{\mbw}   \\
& =  \sum_{I \in \mci} a_I  \mbz^I   =  f(\mbz) = 0
\end{align*}
and  so $\mbw \in \fourIdx{*}{}{}{0}{V}(g)$. 
Thus, 
\[
\left\{ \mbw \in \fourIdx{*}{}{n}{0}{\C} : \std{\mbw} \in V(f)\right\} 
\subseteq \fourIdx{*}{}{}{0}{V}(g). 
\]

Conversely, if $\mbw \in \fourIdx{*}{}{}{0}{V}(g)$, then 
\begin{align*}
0 & = \std{g(\mbw)}  =  \std{\left(\sum_{i \in \mci} b_I \mbw^I\right) } \\
& =  \sum_{i \in \mci} \std{b_I}  \fourIdx{\circ}{}{I}{}{\mbw}    
=  \sum_{i \in \mci} a_I   \fourIdx{\circ}{}{I}{}{\mbw}  \\
& = f( ^{\circ}\mbw )
\end{align*}
and so $^{\circ}\mbw \in V(f)$. 
Thus, 
\[
\fourIdx{*}{}{}{0}{V}(g) \subseteq  \left\{ \mbw \in \fourIdx{*}{}{n}{0}{\C} : \std{\mbw} \in V(f)\right\}. 
\]
This completes the proof. 
\end{proof}

If $n=1$, that is, if $f \in \C[t]$ is a polynomial in one variable, then we recover 
Ross's proof~\cite{ross22} of the continuity of the roots of a polynomial.

\section{Affine varieties}

Let $\mbt = (t_1,\ldots, t_n)$.  
Let $\mcf = \{f_{\lambda}(\mbt) :\lambda \in \Lambda\}$ 
 be a set of polynomials in $\C[\mbt]$. 
The \emph{affine variety}\index{affine variety} defined by \mcf\  is 
\[
V(\mcf)  = \left\{ \mbz \in \C^n: 
f_{\lambda}(\mbz) = 0 \text{ for all } \lambda \in \Lambda \right\}.
\]
Let $\fourIdx{*}{}{}{}{\mcg} = \{g_{\lambda}(\mbt)  :\lambda \in \Lambda\}$ 
 be a set of polynomials in $\fourIdx{*}{}{}{0}{\C} [\mbt]$. 
The \emph{nonstandard affine variety}\index{nonstandard affine variety}\index{affine variety!nonstandard}  defined by $\fourIdx{*}{}{}{}{\mcg}$ is  
\begin{align*}
\fourIdx{*}{}{}{}{V}(\nstd{\mcg}) 
&  = \{\mbw \in \fourIdx{*}{}{n}{}{\C}:  g_{\lambda}(\mbw) \approx 0  \text{ for all } \lambda \in \Lambda\}. 
% &  = \{\mbw \in \fourIdx{*}{}{n}{}{\C}:  \std{g_{\lambda}}(\mbw) = 0   \text{ for all } \lambda \in \Lambda\} \}.
\end{align*}
The \emph{nonstandard finite affine variety} defined by $\fourIdx{*}{}{}{}{\mcg}$ is  
\begin{align*}
\fourIdx{*}{}{}{0}{V}(\nstd{\mcg}) 
&  = \{\mbw \in \fourIdx{*}{}{n}{0}{\C}:  g_{\lambda}(\mbw) \approx 0  \text{ for all } \lambda \in \Lambda\} \\
 &  = \{\mbw \in \fourIdx{*}{}{n}{0}{\C}:  \std{g_{\lambda}}(\mbw) = 0  
 \text{ for all } \lambda \in \Lambda\} \}.
\end{align*}

Let $\fourIdx{*}{}{}{}{\mcg} = \{g_{\lambda}(\mbt)  :\lambda \in \Lambda\}$ 
and $\fourIdx{*}{}{}{}{\mch} = \{h_{\lambda}(\mbt)  :\lambda \in \Lambda\}$
 be  sets of polynomials in $\fourIdx{*}{}{}{0}{\C} [\mbt]$. 
The set  $\fourIdx{*}{}{}{}{\mcg}$ is an infinitesimal deformation of the set $\fourIdx{*}{}{}{}{\mch}$, 
denoted $\fourIdx{*}{}{}{}{\mcg} \approx \fourIdx{*}{}{}{}{\mch}$, if
$g_{\lambda} \approx h_{\lambda}$ for all $\lambda \in \Lambda$. 
Let $\mcf = \{f_{\lambda}(\mbt)  :\lambda \in \Lambda\}$ 
 be a set of polynomials in $\C[\mbt]$.  
 The set  $\fourIdx{*}{}{}{}{\mcg}$ is an infinitesimal deformation of the set $\mcf$, 
denoted $\fourIdx{*}{}{}{}{\mcg} \approx \mcf$, if
$g_{\lambda} \approx f_{\lambda}$ for all $\lambda \in \Lambda$ 
or, equivalently,   if
$\std{g_{\lambda}} = f_{\lambda}$ for all $\lambda \in \Lambda$.

\bt                          \label{NonstandardAffine:theorem:variety}
Let $\mcf = \{f_{\lambda}(\mbt)  :\lambda \in \Lambda\}$ 
 be a set of polynomials in $\C[\mbt]$ and let 
 $\fourIdx{*}{}{}{}{\mcg} = \{g_{\lambda} (\mbt) :\lambda \in \Lambda\}$ 
 be a set of polynomials in $\fourIdx{*}{}{}{0}{\C} [\mbt]$. 
 If  $\fourIdx{*}{}{}{}{\mcg}$ is an infinitesimal deformation of \mcf, 
 then the nonstandard finite affine variety $\fourIdx{*}{}{}{0}{V}(\nstd{\mcg})$ 
 is an infinitesimal deformation of the affine variety $V(\mcf)$. 
 Equivalently, 
 \[
\fourIdx{*}{}{}{0}{V}(\nstd{\mcg}) =  
\left\{ \mbw \in \fourIdx{*}{}{n}{0}{\C} : \std{\mbw} \in V(\mcf)\right\}. 
 \]
\et

\begin{proof}
Every variety is an intersection of hypersurfaces.    
By Theorem~\ref{NonstandardAffine:theorem:hypersurface}, 
\[
\mbw \in \fourIdx{*}{}{}{0}{V}(\nstd{\mcg}) = \bigcap_{\lambda \in \Lambda} \fourIdx{*}{}{}{0}{V}(g_{\lambda})  
\]
if and only if 
\[
\mbw \in  \fourIdx{*}{}{}{0}{V}(g_{\lambda})  \qquad \text{for all $\lambda \in \Lambda$} 
\]
if and only if 
\[
\fourIdx{\circ}{}{}{\lambda}{g}(\mbw) = 0 \qquad \text{for all $\lambda \in \Lambda$} 
\]
if and only if 
\[
f_{\lambda}(\std{\mbw}) = 0 \qquad \text{for all $\lambda \in \Lambda$} 
\]
if and only if 
\[
\std{\mbw} \in V(f_{\lambda})  \qquad \text{for all $\lambda \in \Lambda$} 
\]
if and only if 
\[
\std{\mbw} \in  \bigcap_{\lambda \in \Lambda} V(f_{\lambda}) =  V(\mcf).   
\]
This completes the proof. 
\end{proof}

\section{Open problems} 
The distance from a point \mbw\ to a nonempty set $Z$ in a metric space is 
\[
d(\mbw,Z) =  \inf \{ d(\mbw, \mbz) : \mbz \in Z\}. 
\]
The \emph{Hausdorff distance}\index{Hausdorff distance} between nonempty 
sets $W$ and $Z$ is 
\[
d_H(W,Z) = \max\left( \sup_{\mbw\in W} d(\mbw,Z), \sup_{\mbz \in Z} d(\mbz, W)  \right).  
\]

For $\mbz = (z_1,\ldots, z_n) \in \C^n$, let $\|\mbz\|_{\infty} = \max \left\{|z_i|: i \in \{1,\ldots, n\} \right\}$. 
The set $\C^n$ is a  metric space with distance function $d(\mbw,\mbz) = \|\mbw-\mbz\|_{\infty}$. 
Let $W$ and $Z$ be nonempty subsets of $\C^n$.  
The set $W$ is an \emph{$\varepsilon$-deformation}\index{deformation!$\varepsilon$-deformation}  
of $Z$, denoted $W \approx_{\varepsilon} Z$, if $d_H(W,Z) < \varepsilon$.    
Equivalently, $W$ is an $\varepsilon$-deformation of $Z$ 
if for each $\mbw \in W$ there exists $\mbz \in Z$ with $\|\mbw - \mbz\|_{\infty} < \varepsilon$   
and for each $\mbz \in Z$  there exists $\mbw \in W$ with $\|\mbw - \mbz\|_{\infty} < \varepsilon$. 
This relation is symmetric:  $W \approx_{\varepsilon} Z$ if and only if $Z \approx_{\varepsilon} W$. 

Let $\mbt = (t_1,t_2)$ and $0 < \delta' < \delta$.  Consider the polynomials 
\[
f(\mbt) = t_2-t_1 
\]
and 
\[
g(\mbt) = t_2 - (1+\delta')t_1  
\]
in the ring $\C[\mbt] = \C[t_1,t_2]$.  
The polynomial $g(\mbt)$ is a $\delta$-deformation of the polynomal $f$. 
We obtain the hypersurfaces 
\begin{align*}
V(f) & = \{\mbz = (z_1,z_2) \in \C^2: f(\mbz) = 0 \} \\ 
& = \{\mbz = (z,z):z\in \C\} 
\end{align*}
and 
\begin{align*}
V(g) & = \{\mbw = (w_1,w_2) \in \C^2: g(\mbw) = 0 \}  \\
& = \{\mbw = (w,(1+\delta')w):w\in \C\}. 
\end{align*}
For all $\mbz  = (z,z) \in V(f)$ and $\mbw = (w,(1+\delta')w) \in V(g)$ we have 
\[
\|\mbw - \mbz\|_{\infty} = \max\left(|w-z|, |w-z+\delta' w | \right). 
\]

Let $\varepsilon > 0$.  
If $\|\mbw - \mbz\|_{\infty}  < \varepsilon$, then 
\[
|w-z| < \varepsilon \qqand |w-z + \delta' w | < \varepsilon
\]
and so  
\[
|  \delta' w | - |w-z| \leq |w-z + \delta' w | < \varepsilon 
\]
and 
\[
|  \delta' w |  < \varepsilon +  |w-z|  < 2\varepsilon. 
\]
It follows that 
\[
|w| < \frac{2\varepsilon}{\delta'}.  
\]
Equivalently, for all $\varepsilon > 0$ and $0 < \delta' < \delta$, if $w \in \C$ and 
\[
|w| \geq \frac{2\varepsilon}{\delta'}   
\]
then 
\[
\|\mbw - \mbz\|_{\infty}  \geq \varepsilon 
\]
for all $\mbz \in V(f)$.  
Thus, if $\varepsilon > 0$, then for every $\delta > 0$ the polynomial $g(\mbt)$ 
is a $\delta$-deformation of the polynomial $f(\mbt)$ but the hypersurface $V(g)$ 
is not an $\varepsilon$-deformation of the hypersurface $V(f)$.

It would be of interest to translate Theorems~\ref{NonstandardAffine:theorem:hypersurface} 
and~\ref{NonstandardAffine:theorem:variety} from the language of the nonstandard 
mathematics to standard mathematics.  The simple example above shows that 
it is not true that, for every $\varepsilon > 0$, 
there exists $\delta > 0$ such that 
if the polynomial $g \in \C[\mbt]$ is a $\delta$-deformation of  the polynomial 
$f \in \C[\mbt]$, then the hypersurface 
$V(g)$ is an $\varepsilon$-deformation of the hypersurface $V(f)$.

There is the following elementary result. 
For $I = (i_1,\ldots, i_n) \in \N_0^n$, let $|I| = \sum_{j=1}^n i_j$. 
Let $|\mci|$ denote the cardinality of the finite set \mci. 
For all $T > 0$ we define the $n$-dimensional complex \emph{hypercube}
\[
H(T) = \{ \mbz \in \C^n: \| \mbz \|_{\infty} \leq T \}.
\]

\bl                                    \label{NonstandardAffine:lemma:standard-epsilon}
Let $\mbt = (t_1,\ldots, t_n)$, let $\mci$ be a nonempty finite subset of $\N_0^n$,  and let 
\[
f(\mbt) = \sum_{I \in \mci} a_I \mbt^I 
\] 
be a nonzero polynomial in $\C[\mbt]$.  
Let 
\[
d = \max( |I|:I \in \mci).
\]
Let $T \geq 1$.  For every $\varepsilon > 0$, let 
\[
0 < \delta < \frac{\varepsilon}{T^d |\mci|}.
\]
Let the polynomial 
\[
g(\mbt) = \sum_{I \in \mci} b_I \mbt^I \in \C[\mbt]
\]
be a $\delta$-deformation of $f$.  
If $\mbz \in H(T)$, then 
\[
|g(\mbz)-f(\mbz)| < \varepsilon. 
\]
\el

\begin{proof}
Let $\mbz = (z_1,\ldots, z_n) \in \C^n$ with 
\[
\|\mbz\|_{\infty} = \max(|z_i|:i=1,\ldots, n ) \leq T.
\] 
For all $I = (i_1,\ldots, i_n)  \in \mci$ we have 
\[
|b_I - a_I| < \delta
\]
and 
\[
\left|  \mbz^I   \right| = \left|  z_1^{i_1} \cdots  z_n^{i_n}    \right| 
\leq    T^{i_1+ \cdots  +i_n}  =  T^{|I|} \leq T^d. 
\]
It follows that 
\begin{align*}
|g(\mbz)-f(\mbz)| 
& = \left|  \sum_{I \in \mci} (b_I - a_I) \mbz^I   \right| \\
& \leq  \sum_{I \in \mci} |b_I - a_I|  \left|  \mbz^I   \right| \\
& \leq \delta T^d  |\mci|  \\ 
& < \varepsilon.
\end{align*} 
This completes the proof.  
\end{proof}

Let $g(\mbt)$ be a nonzero polynomial in $\C[\mbt]$.  
For all $\varepsilon > 0$, we define  
\[
V_{\varepsilon}(g) = \{\mbw \in \C^n: |g(\mbw)| < \varepsilon \}. 
\]

\bt                                     \label{NonstandardAffine:theorem:standard-epsilon}
Let $\mbt = (t_1,\ldots, t_n)$, let $\mci$ be a finite subset of $\N_0^n$,  and let 
\[
f(\mbt) = \sum_{I \in \mci} a_I \mbt^I 
\] 
be a nonzero polynomial in $\C[\mbt]$.  For all $\varepsilon > 0$ and $T > 0$ 
there exists $\delta > 0$ such that if the polynomial $g(\mbt) \in \C[\mbt]$ 
is a $\delta$-deformation of $f(\mbt)$, then 
\[
V(f) \cap H(T)  \subseteq V_{\varepsilon}(g). 
\]
\et

\begin{proof}
This follows from Lemma~\ref{NonstandardAffine:lemma:standard-epsilon}.  
If $\mbz \in V(f) \cap H(T)$, then $f(\mbz) = 0$ and so 
$|g(\mbz)|  = |g(\mbz)-f(\mbz)| < \varepsilon$ and $\mbz \in V_{\varepsilon}(g)$. 
This completes the proof. 
\end{proof}

Another interesting open problem is to investigate the stability of multiplicities 
of points on hypersurfaces and affine varieties under infinitesimal deformations.

\section{Acknowledgements}
This paper was inspired by several  conversations with David Ross 
and by his recent paper~\cite{ross22}.  
I also thank Renling Jin for many discussions about nonstandard analysis over many years 
and for his paper~\cite{jin08}.

\end{document}